\documentclass[a4paper, 12pt]{article}
\usepackage{enumerate, theorem}
\usepackage{amsmath, amsfonts, amssymb}
\usepackage[height=22.5cm, width=15cm]{geometry}
\usepackage[all]{xy}
\usepackage{mathrsfs, yfonts}

\newcommand{\K}{\hat{\mathcal{K}}}
\newcommand{\h}{\hat{\mathcal{H}}}

\newtheorem{thm}{Th\'eor\`eme}[subsection]
\newtheorem{defn}[thm]{D\'efinition}
\newtheorem{cor}[thm]{Corollaire}
\newtheorem{prop}[thm]{Proposition}
\newtheorem{lemma}[thm]{Lemme}

\begin{document}


\title{On the Brieskorn (a,b)-module of an isolated hypersurface singularity.} 

 \author{Daniel Barlet.}
 
 \date{}
 
 \maketitle
 
 \markright{On the Brieskorn (a,b)-module ...}
 
 \section*{Abstract}
 
 We show in this note that for a germ \ $g$ \  of holomorphic function with an isolated singularity at the origin of \ $\mathbb{C}^n$ \ there is a pole for the meromorphic extension of the distribution
 \begin{equation*}
   \frac{1}{\Gamma(\lambda)} \int_X \vert g \vert^{2\lambda}\bar{g}^{-n} \square  \tag{*}
  \end{equation*}
  at \ $- n - \alpha$ when \ $ \alpha$ \ is the smallest root in its class modulo \ $\mathbb{Z}$ \ of the reduce Bernstein-Sato polynomial of \ $g$. This is rather unexpected result comes from the fact that the self-duality of the Brieskorn (a,b)-module \ $E_g$ \ associated to \ $g$ \ exchanges the biggest simple pole sub-(a,b)-module of \ $E_g$ \ with the saturation of \ $E_g$ \ by \ $b^{-1}a$.
  
\noindent In the first part of this note, we prove that the biggest simple pole sub-(a,b)-module of the Briekorn (a,b)-module \ $E$ \ of \ $g$ \ is "geometric" in the sense that it depends only on the hypersurface germ \ $\{ g = 0 \}$ \ at the origin in \ $\mathbb{C}^n$ \ and not on the precise choice of the reduced equation \ $g$, as the poles of (*).\\
 By duality, we deduce the same property for the saturation \ $\tilde{E}$ \ of \ $E$. This duality gives also the relation between the "dual" Bernstein-Sato polynomial and the usual one, which is the key of the proof of the theorem.

  \bigskip
  
  \noindent {\bf Key words} Isolated hypersurface singularity, Brieskorn (a,b)-module, Bernstein-Sato polynomial, dual Bernstein-Sato polynomial.
  
  \smallskip
  
  AMS Classification : 32-S-05, 32-S-25, 32-S-40.
  
  \newpage
 
 \section{Introduction.}
 
 \noindent Let \ $\tilde{g} : (\mathbb{C}^n, 0)  \to (\mathbb{C}, 0) $ \ a germ of holomorphic function with an isolated singularity. Denote by \ $ g : X \to D $ \ a Milnor representative of \ $\tilde{g}$.\\
 Let \ $b_g$ \ be the reduced Bernstein-Sato polynomial of \ $g$. Let \ $\alpha$ \ be the biggest root of \ $b_g$ \ in its class modulo \ $\mathbb{Z}$.  A classical question is whether for \ $j \in \mathbb{N}$ \ big enough the meromorphic extension of the distribution
 $$\frac{1}{\Gamma(\lambda)} \int_X \vert g \vert^{2\lambda}\bar{g}^{-j} \square  $$
 has a pole at \ $\lambda = \alpha $.\\
 The present note gives a result which, in a sense, suggests that, may be, this question is not the good one.\\
  Let me  introduce the dual Bernstein-Sato polynomial of \ $g$ \ by the formula
 $$ b^*_g(z) = (-1)^q.b_g(- n - z) $$
 where \ $q : = deg(b_g)$. Recall that all roots of \ $b_g$ \ (and $b_g^*$) \ are contained in \ $]-n, 0[ $, see [K.76] for the inequality \ $< 0$, and the section 3 for the inequality \ $> -n$.
 
 \smallskip
 
\noindent  We shall prove the following result.
 
 \begin{thm}\label{Cont.pr.}
 Let \ $\alpha$ \ be the \underline{smallest}\footnote{recall that we are dealing with negative numbers.} root of \ $b_g$ \ in its class modulo \ $\mathbb{Z}$, and let \ $d$ \ be its multiplicity (as a root of \ $b_g$). Then the meromorphic extension of the of the distribution
 $$\frac{1}{\Gamma(\lambda)} \int_X \vert g \vert^{2\lambda}\bar{g}^{-n} \square  $$
 has a p\^ole of order \ $\geq d$ \ at \ $-n - \alpha$.
 \end{thm}
 
 \bigskip
 
 \noindent {\bf Remarks.}
 \begin{enumerate}
 \item In general \ $b^*_g \not= b_g $ \ so it is not clear that \ $- n -\alpha$ \ is a root of \ $b_g$. But, of course, the previous theorem implies that there exists at least \ $d$ \ roots of \ $b_g$ (counting multiplicities) which are bigger than \ $-n - \alpha$. If \ $- n - \alpha \in [-1,0[$ \ then there is no choice : \ $- n - \alpha$ \ is a root of multiplicity $\geq d$ \ of \ $b_g$.
 \item This result gives, in term of the Bernstein-Sato polynomial \ $b_g$, a precise value where we know that a pole appears in the class \ $[\beta]$ \ modulo \ $\mathbb{Z}$ \ of a root \ $\beta$ \ of \ $b_g$. But the pole which is given is not at the biggest root of \ $b_g$ \ in this class but at the biggest root of \ $b^*_g$ \ in this class ! \\
 A clear reason for that is given in the proof: the dual Bernstein-Sato polynomial is the minimal polynomial of \ $-b^{-1}a$ \ acting on \ $F/b.F$ \ where \ $F$ \ is the biggest simple pole sub-(a,b)-module of the Brieskorn (a,b)-module \ $E$ \ associated to \ $g$. So it lies in the lattice given by \underline{holomorphic} forms.\\
  On the contrary, \ $b_g$ \ is the  the minimal polynomial of \ $-b^{-1}a$ \ acting on \ $\tilde{E}/b\tilde{E}$ \ where \ $\tilde{E}$ \ is the saturation of \ $E$ \ by \ $b^{-1}a$, or, in other words, the minimal simple pole (a,b)-module containing \ $E$. So, if \ $E$ \ is not a simple pole (a,b)-module, elements in \ $\tilde{E}$ \ are not always representable in the holomorphic lattice, and so we may need some power of \ $g$ \ as denominators. And this may introduce integral shifts for the poles.
 \item The case where \ $E$ \ is a simple pole (a,b)-module (that is to say  when we have \ $F = E = \tilde{E}$) corresponds to a quasi-homogeneous \ $g$, with a suitable choice of coordinates. In this case we have \ $b_g^* = b_g$, so \ $- n - \alpha$ \ is the smallest root of \ $b_g$ \ in its class modulo \ $\mathbb{Z}$.
\end{enumerate}
 
 \noindent In the first part of this note, we prove that the biggest simple pole sub-(a,b)-module of the Briekorn (a,b)-module \ $E$ \ of \ $g$ \ is "geometric" in the sense that it depends only on the hypersurface germ \ $\{ g = 0 \}$ \ at the origin in \ $\mathbb{C}^n$ \ and not on the precise choice of the reduced equation \ $g$.\\
Remark that the poles of the meromorphic distributions \ $\frac{1}{\Gamma(\lambda)} \int_X \vert g \vert^{2\lambda}\bar{g}^{-j} \square  $ are also "geometric" in the sense above.\\
By duality, we deduce the same property for the saturation \ $\tilde{E}$ \ of \ $E$. This duality gives also the relation between the dual Bernstein-Sato polynomial and the usual one, which is the key of the proof of the theorem.

 \section{Changing the reduced equation.}
 
 \noindent Let \ $ g : X \to D$ \ be a Milnor representative of a germ of an holomorphic function with an isolated singularity at the origin of  \ $\mathbb{C}^n, n \geq 2 $. We define the function
  $$ f(t, x) : = e^t.g(x) \quad {where} \quad  f : \mathbb{C} \times X \to \mathbb{C} $$
  and we denote by  \ $\pi : \mathbb{C} \times \mathbb{C} \times X$ \ the projection defined by \ $ \pi(\lambda, t, x) = (t, x) $. We shall denote by \ $F$ \ the function \ $\pi^*(f)$. Its critical locus is \ $S : = \mathbb{C} \times \mathbb{C} \times \{0\} $.\\
We consider on \ $Y = \{ F = 0\}$, as in [B.05], the complex of sheaves   $\big((\hat{\mathcal{K}}er\,dF)^{\bullet}, d^{\bullet}\big)$. The following theorem is an easy generalization of [B.05] th.2.2 (case LII) .
  
  \begin{thm}\label{loc.cst.}
 In the situation describe above, the $n-$th cohomology sheaf of the complex $\big((\hat{\mathcal{K}}er\,dF)^{\bullet}, d^{\bullet}\big)$ \ is a constant sheaf whose fiber is $F_g$  the biggest simple pole sub-(a,b)-module of the Brieskorn (a,b)-module $E_g$ associated to the function $g$.
  \end{thm}
 
 \bigskip
 
 \noindent It is easy to deduce from the previous theorem the following corollary.
 
 \begin{cor}
 Let \ $g$ \ be a germ of an holomorphic function with an isolated singularity at the origin of \ $\mathbb{C}^n$. Let \ $h$ \ be  any invertible holomorphic  germ  at the origin. Then  the biggest  simple pole sub-(a,b)-module of the Brieskorn (a,b)-module associated to the function \ $h.g$ \ does not depend on the choice of \ $h$ \ up to isomorphism.\\
 More precisely, if the holomorphic invertible function depends holomorphically on some parameter \ $\lambda$ \  in a complex manifold \ $\Lambda$, the subsheaf of the sheaf on \ $\Lambda$ \ defined by the Brieskorn (a,b)-modules of the fibers\footnote{we defini this sheaf via  the  cohomology of the formal completion of the de Rham complex of $\Lambda-$relative holomorphic forms annihilated by \ $\wedge dF$.}, which is given in each fiber by  the biggest  simple pole sub-(a,b)-module of the Brieskorn (a,b)-module, is a locally constant sheaf on \ $\Lambda$.
 \end{cor}
 
 \bigskip
 
 \noindent \textit{Proof of the theorem.} Let us first consider the case of an holomorphic function \ $f$ \ on a complex manifold \ $Z$ \ and let  the holomorphic function \ $F$ \ be  \ $F : = \pi^*(f)$ \ on \ $\mathbb{C} \times Z$ \  where \ $\pi :  \mathbb{C} \times Z \to Z$ \ is the projection.\\
 In this situation we have the following description of \ $(\K er\, dF)^p$ :
 $$ (\K er\, dF)^p = \pi^*((\K er\,df)^p) \oplus d\lambda\wedge \pi^*((\K er\,df)^{p-1}).$$
 Then \ $ \alpha \oplus d\lambda\wedge\beta \in (\K er\, dF)^p$ \ is $d-$closed iff it satisfies :
 $$d_/ \alpha = 0 \quad {\rm and} \quad \frac{\partial \alpha}{\partial \lambda} = d_/ \beta $$
 where \ $\frac{\partial \alpha}{\partial \lambda}$ \ is defined by the equation \ $d\alpha = d_/\alpha + d\lambda \wedge \frac{\partial \alpha}{\partial \lambda}$.\\
 
 \begin{lemma}
 In the situation above set \ $Y = \{ f = 0\}$ ;  we have the short exact sequence of complex of sheaves on \ $\mathbb{C} \times Y $: 
   $$ 0 \to (\K er\, dF^{\bullet}, d^{\bullet}) \to \big(\pi^*(\K er\, df^{\bullet}), d_/^{\bullet}\big) \overset{\frac{\partial}{\partial \lambda}}{\to} \big(\pi^*(\K er\, df^{\bullet}), d_/^{\bullet}\big) \to 0 .$$
  So if the sheaf \ $\h^{p-1}_f$ \ is \  $0$ \ on \ $Z$ \ for \ $p \geq 3$ \ or is isomorphic to\footnote{recall that \ $E_1 : = \mathbb{C}[[b]].e_1$ \ where \ $a.e_1 = b.e_1$.} \ $E_1 \otimes \mathbb{C}_Y$ \  for \ $p = 2$, then we have for \ $p \geq 2$ \  the exact sequence of sheaves on  \ $\mathbb{C} \times Y$ :
 $$ 0 \to \h^p_F \overset{i}{\to}  \pi^*(\h^p_f) \overset{\partial/\partial \lambda}{\longrightarrow} \pi^*(\h^p_f) .$$
 \end{lemma}

\bigskip

 \noindent \textit{Proof.} Here the sheaf \ $\pi^*(\h^p_f)$ \ is defined via $\lambda-$relative holomorphic forms. On this complex we have a derivation  \ $\partial/\partial \lambda$ \ commuting with the product by the function \ $F$, the wedge product with \ $dF$ \ and the \ $\lambda-$relative de Rham differential denoted by  \ $d_/$. Remark also that we have \ $d_/F = dF $. \\
 The exactness of the short exact sequence of complexes is obvious and the associated  long exact  cohomology sequence is enough to conclude for \ $p \geq 3$. For the \ $p = 2$ \ case, we have only to check the injectivity of the map \ $i $.\\
  Let \ $ \alpha \oplus d\lambda \wedge \beta \in  (\K er\, dF)^p \cap Ker\, d $; its image by \ $i$ \ is the class \ $ [\alpha]$. If it vanishes in  \ $\pi^*(\h^p_f)$ \ we can find \ $\gamma \in \pi^*((\K er\, df)^{p-1})$ \ such that\ $d_/\gamma = \alpha$. Differentiating with respect to \ $\lambda$ \ gives, using the relation \ $\frac{\partial \alpha}{\partial \lambda} = d_/ \beta $,
  $$d_/(\beta - \frac{\partial \gamma}{\partial \lambda}) = 0 .$$
  But as \ $\beta - \frac{\partial \gamma}{\partial \lambda} \in \pi^*((Ker\, df)^{p-1})$ \ this form induces a class in \  $\pi^*(\h^{p-1}_f)$. So we can write 
   $$ \beta = \frac{\partial \gamma}{\partial \lambda} + \varphi(\lambda,f).df $$
   where \ $\varphi \in \pi^*(\mathbb{C}[[z]])$. We obtain, if \ $\frac{\partial \psi}{\partial \lambda}(\lambda, f) = \varphi(\lambda, f)$ :
    $$\alpha + d\lambda \wedge \beta = d(\gamma +  \psi(\lambda,f).df)$$
    which allows to conclude, as \ $\gamma +  \psi(\lambda,f).df$ \ is in  \ $  \pi^*((\K er\, df)^{1}). \hfill \blacksquare$

\bigskip
 
  \noindent \textit{End of the proof of the theorem.} We proved in [B-05] theorem 2.2 that  the sheaf \ $\h^n_f$ \ is a constant sheaf on \ $\mathbb{C}\times \{0\} \subset \mathbb{C}\times X = Z$ \ with fiber the biggest simple pole sub-(a,b)-module in \ $E_g$. So the sams is true for the sheaf \ $\h^n_F$ \ on \ $\mathbb{C}\times\mathbb{C}\times \{0\}. \hfill \blacksquare$

 \bigskip
 
 \noindent \textit{Proof of the corollary.} Let \ $ c : \mathbb{C}\times X \to \mathbb{C}$ \ be an holomorphic function and set \ $h_{\lambda}(x) : = exp(c(\lambda, x))$ \ for \ $ \lambda \in \mathbb{C}$ \ and \ $x \in X$. Choose the following coordinate system on $\mathbb{C}\times\mathbb{C}\times X$ \ near the point \ $(\lambda_0, t_0, 0)$ :
  $$ \lambda' = \lambda, \quad t' = t - c(\lambda, x), \quad x' = x .$$
 Then the function \ $F$ \ is  transformed in \ $\tilde{F}(\lambda', t', x') = e^{t'}.(e^{c(\lambda', x')}.g(x')) = F(\lambda, t, x)$. The corollary follows, because we can always join two invertible functions inside an holomorphic family of invertible functions (and the restriction of a constant sheaf is a constant sheaf).$\hfill \blacksquare$

  \bigskip
 
 \section{The dual Bernstein-Sato polynomial.}
 
 \noindent We shall now consider an (a,b)-module \ $E$ \ such that
  \begin{enumerate}[i)]
 \item The (a,b)-module \ $E$ \ is regular (see [B.93]).
\item There exists a complex number \ $\delta$ \ and an isomorphism of (a,b)-modules \\ 
$ \kappa : \check{E} \to Hom_{a,b}(E, E_{\delta})$, where \ $\check{E}$ \ is the (a,b)-module \ $E$ \ in which "$a$" and "$b$" acts as \ $- a$ \ and \ $- b$.
\end{enumerate}

\noindent Recall, for the convenience of the reader, that \ $E_{\delta}$ \ is the rank 1 (a,b)-module defined by  \ $E_{\delta} : = \mathbb{C}[[b]].e_{\delta}$ \ where $a$ \ acts by \ $a.e_{\delta} = \delta.b.e_{\delta}$. \\
If \ $E$ \ and \ $F$ \ are (a,b)-modules, the (a,b)-module \ $Hom_{a,b}(E, F)$ \ is defined as follows : we define on the \ $\mathbb{C}[[b]]-$module  \ $Hom_{\mathbb{C}[[b]]}(E, F)$,  which is free and of finite rank, an action of \ $a$ \ by the formula :
$$ (a.\varphi)(x) =  a_F.\varphi(x) -  \varphi(a_E.x), \quad \forall x \in E .$$
Of course, we have to check that \ $a.\varphi$, defined in this way, is  \ $\mathbb{C}[[b]]-$linear and that we have  \ $a.b.\varphi - b.a.\varphi = b^2.\varphi $. It is not difficult to check also that \ $Hom_{a,b}(E, F)$ \ is regular when \ $E$ \ and \ $F$ \ are regular (see [B.95]).

\noindent Recall also that the Brieskorn (a,b)-module of a germ of holomorphic function with an isolated singularity in \ $\mathbb{C}^{n}$ \ satisfies properties i) and ii) above with \ $ \delta = n $, see [Be.01].

\begin{prop}\label{Dual}
Under hypotheses i) and ii) above, let \ $F$ \ be the biggest simple pole sub-(a,b)-module in \ $E$, and let \ $\tilde{E}$ \ the saturation of \ $E$ \ for \ $b^{-1}a$.\\
Then we have natural isomorphisms of (a,b)-modules deduced from \ $\kappa$ :
$$\kappa' :\check{ \tilde{E}} \to Hom_{a,b}(F, E_{\delta}) \quad {\rm and} \quad \kappa'' : \check{F} \to  Hom_{a,b}(\tilde{E}, E_{\delta}) .$$
\end{prop}

\bigskip

\noindent In the proof of this proposition we shall use the following  lemmas.

\begin{lemma}\label{s.p.}
Let \ $E$ \ and \ $F$ \ be simple pole (a,b)-modules. Then  \ $Hom_{a,b}(E, F)$ \  is also a simple pole (a,b)-module.
\end{lemma}

\bigskip

\noindent\textit{Proof.} Fix an element \ $\varphi \in Hom_{a,b}(E, F)$. Then define   \ $\theta : E \to F $ \ by the formula \ $\theta(x) : =  b^{-1}.a.\varphi(x) -  b^{-1}.\varphi(a.x)$ \ for all  \ $x \in E$. As \ $E$ \ has a simple pole, we have \ $a.x \in b.E$ \ and so \ $\varphi(a.x) \in b.F $ \ from $b-$linearity of \ $\varphi$. But \ $F$ \ has also a simple pole, so \ $ b^{-1}.a : F \to F $ \  is well defined.\\
Now \ $\theta$ \ is $b-$linear :
\begin{align*}
\theta(b.y) & =  b^{-1}.a.\varphi(b.y) - b^{-1}.\varphi(a.b.y) = (a + b).\varphi(y) - \varphi((a + b).y) \\
\quad         & =  a.\varphi(y) - \varphi(a.y) = b.\theta(y).
\end{align*}
But we have \ $a.\varphi = b.\theta$ \ in \ $Hom_{a,b}(E, F)$. Therefore \ $Hom_{a,b}(E, F)$ \ is a simple pole (a,b)-module. $\hfill \blacksquare$

\bigskip

\begin{lemma}\label{bidual}
Let \ $E$ \ be a regular (a,b)-module and let  \ $\delta $\ be any  complex number. Then we have a canonical (a,b)-module isomorphism
$$\tau :  E \to Hom_{a,b}( Hom_{a,b}(E, E_{\delta}), E_{\delta}).$$
\end{lemma}

\noindent \textit{Proof.} The map \ $\tau$ \ is defined by \ $x \to \tau(x)[\varphi] = \varphi(x) $. It is obviously a \ $b-$linear isomorphism. So we have only to check the $a-$linearity. But, with the notation \ $\theta = \tau(x)$, we have :
$$(a.\theta)[\varphi] =  a.\big(\theta[\varphi]\big)- \theta[a.\varphi] =a.\varphi(x) - \big(a.\varphi(x) - \varphi(a.x)\big) = \tau(a.x)[\varphi] .$$
And so  \ $a.\tau(x) = \tau(a.x)$. $\hfill \blacksquare$

\bigskip

\begin{lemma}\label{check}
Let \ $E$ \ and \ $F$ \ be two (a,b)-modules.  Then we have a canonical isomorphism
$$Hom_{a,b}(E, F)\,\check{} \to Hom_{a,b}(\check{E}, \check{F}).$$
\end{lemma}

\bigskip

\noindent\textit{Proof.} It is clear that \ $Hom_{a,b}(\check{E}, \check{F})$ \ is the same complexe vector space than \ $Hom_{a,b}(E, F)$ \ and that the action of $b$ \ on it is given by \ $-b$. The fact that the action of $a$ \ is the opposite of the action of $a$ on $Hom_{a,b}(E, F)$ \ follows also directly from the definition of \ $Hom_{a,b}. \hfill \blacksquare$

\bigskip

\noindent\textit{Proof of proposition \ref{Dual}.} The functor \ $Hom_{a,b}(- , E_{\delta})$ \ applied to the inclusion of \ $E$ \ in \ $\tilde{E}$ \ gives an (a,b)-linear injection
$$ Hom_{a,b}(\tilde{E}, E_{\delta}) \hookrightarrow Hom_{a,b}(E, E_{\delta})\simeq  \check{E} .$$
As \ $ Hom_{a,b}(\tilde{E}, E_{\delta})$ \ has a simple pole by lemma \ref{s.p.} it is contained in \ $\check{F}$, by definition of \ $F$. Apply now the functor \ $Hom_{a,b}(-,E_{\delta})$ \ to the inclusions 
$$  Hom_{a,b}(\tilde{E}, E_{\delta}) \hookrightarrow \check{F} \hookrightarrow \check{E} $$

This gives (a,b)-linear injections
$$  Hom_{a,b}(\check{E}, E_{\delta}) \hookrightarrow Hom_{a,b}( \check{F},  E_{\delta}) \hookrightarrow \tilde{E} $$
using lemma \ref{bidual}. But, as \ $\check{E}_{\delta}$ \ is canonically isomorphic to \ $E_{\delta}$, so  we have isomorphims
$$  Hom_{a,b}(\check{E}, E_{\delta}) \simeq  Hom_{a,b}(\check{E}, \check{E}_{\delta})\simeq Hom_{a,b}(E, E_{\delta})\,\check{} \simeq \check{\check{E}} \simeq E $$
using lemma \ref{check} and our hypothesis on \ $E$. So the simple pole (a,b)-module \ $Hom_{a,b}( \check{F},  E_{\delta})$ \ which lies between \ $E$ \ and \ $\tilde{E}$ \ is equal to  \ $\tilde{E}$. We conclude using again the canonical isomorphism between \ $E_{\delta}$ \ and \ $\check{E}_{\delta}$ \ and the lemma \ref{bidual}. $\hfill \blacksquare$

\bigskip

\noindent{\bf Remark.}

\smallskip 

\noindent In the situation of the proposition  \ref{Dual} the non-degenerate (a,b)-bilinear pairing
$$ h : \check{E} \times E \to E_{\delta}$$
deduced from \ $\kappa$ \ via the formula \ $h(x,y) : = \kappa(x)[y] $, gives also  non-degenerate (a,b)-bilinear pairings 
$$ h' : \check{\tilde{E}} \times F \to E_{\delta} \quad {\rm and} \quad h'' : \check{F} \times \tilde{E} \to E_{\delta} $$
deduced from \ $\kappa'$ \ and \ $\kappa''$ \ via the formulas \ $h'(x,y) :  = \kappa'(x)[y]  $ \ and \ $h''(u,v) = \kappa''(u)[v]$.

\bigskip

An obvious consequence of proposition \ref{Dual}  is the following corollary of the theorem \ref{loc.cst.}.

\begin{cor}
Let \ $g$ \ be a germ of an holomorphic function having an isolated singularity at the origin in \ $\mathbb{C}^n$ \ where \ $n \geq 2$. For any holomorphic invertible germ \ $h$ \ at the origin, the saturation by \ $b^{-1}a$ \ of the Brieskorn (a,b)-module of the germ \ $h.g$ \ is independant, up to an isomorphism of (a,b)-module, of the choice of \ $h$.\\
If the invertible  \ $h$ \  depends holomorphically of a parameter  \ $\lambda$ \ in a complex manifold  \ $ \Lambda$, the sheaf on  \ $ \Lambda$ \ defined by the saturations of the Brieskorn (a,b)-modules of the germs \ $h_{\lambda}.g$ \ is a locally constant sheaf on \ $\Lambda$. 
\end{cor}

\bigskip

\section{Poles of \ $\int_X \vert g \vert ^{2.\lambda} \square .$}

\noindent We shall begin by a simple definition.

\bigskip

\begin{defn}
Let \ $E$ \ be a regular (a,b)-module. We shall call {\bf dual Bernstein polynomial} of \ $E$, denoted by  \ $b^*_E$, the minimal polynomial of the linear endomorphism  \ $- b^{-1}.a$ \ acting on the (finite dimensional) vector space \ $F / b.F$ \ where \ $F$ \ is the biggest simple pole sub-(a,b)-module of \ $E$.
 \end{defn}
 
 \bigskip
 
 \noindent Recall that the Bernstein-Sato polynimial of \ $E$ \ is the minimal polynomial of the action of \ $- b^{-1}.a$ \ on the (finite dimensional) vector space  \ $\tilde{E}/b.\tilde{E}$, where \ $\tilde{E}$, as before, is the saturation of \ $E$ \ by \ $b^{-1}.a$. In other words, \ $\tilde{E}$ \ is the smallest simple pole (a,b)-module which contains \ $E$. This can be understood in two ways. Either you look in \ $E[b^{-1}]$ \ for the smallest simple pole (a,b)-module containing \ $E$. The other way is to consider the inclusion \ $E \to \tilde{E}$ \ as the initial element for inclusions of \ $E$ \ in simple poles (a,b)-modules.
 
  \bigskip
 
 \noindent {\bf Remark.}
 
 \smallskip
 
 \noindent Let \ $\delta $ \ a given complex number, and assume that the (a,b)-module  \ $E$ \ is equipped with an (a,b)-linear  isomorphism
 $$ \kappa : \check{E} \to Hom_{a,b}(E, E_{\delta}) .$$
Then we have  \ $b^*_E(z) = (-1)^r.b_E(- \delta - z) $ \ where \ $ r : = deg(b_E)$, since \ $b^{-1}a$ \ acts on the same way on \ $E$ \ and \ $\check{E}$.\\
So, for the Brieskorn (a,b)-module of a germ of an holomorphic function \ $g$ \  with an isolated singularity at the origin of  \ $\mathbb{C}^n$ \ the dual Bernstein polynomial is given by
$$ b^*_g(z) = (-1)^{r}b_g(- n - z) .$$
Using Malgrange positivity theorem it is easy to show that the roots of \ $b_g^*$ \ are strictly negative. This gives, using [K.76],  the fact that the roots of \ $b_g$ \ are contained in \  $]- n, 0[$.

\bigskip

\noindent \textit{Proof of the theorem \ref{Cont.pr.}} The only new point for this proof, compared to [B.84 a] and  [B.84 b], is the following :\\
In a simple pole (a,b)-module \ $F$,  if a spectral value \ $\beta$ \ of multiplicity \ $d$ \ for the action of \ $b^{-1}.a$ \ on \ $F/bF$, is minimal in its class modulo \ $\mathbb{Z}$, there exists elements \ $e_1, \cdots, e_d$ \ in \ $F$, giving a Jordan block of size \ $d$ \ for \ $b^{-1}a$ \ acting on \ $F/bF$, and such that they satisfy in \ $F$ \ the relations
$$ a.e_j = \beta.b.e_j + b.e_{j-1}, \ \forall j \in [1,d]$$
with the convention \ $e_0 = 0$ (see [B.93]).\\
This enable us, using the standard technics of [B.84 a], to build up \ $(n-1)-$holomorphic forms \ $\omega_1, \cdots, \omega_d$ \  in a neighbourghood of the origin in \ $\mathbb{C}^n$, such that
$$ d\omega_j = \beta.\frac{dg}{g}\wedge \omega_j +  \frac{dg}{g}\wedge \omega_{j-1}, \forall  j \in [1,d]$$
with the convention \ $\omega_0 = 0 $, which induce a Jordan block of size \ $d$ \ in \ $H^{n-1}(F, \mathbb{C})$ \ where \ $F$ \ is the Milnor fiber of  \ $g$, for the eigenvalue  \ $exp(-2i\pi.\beta)$ \ of the monodromy.\\
So we avoid in this way the integral  shifts coming from the use of  a lattice which may be not contained in the one given by holomorphic forms and we can realize the pole of our statement for \ $\lambda = - \beta$, using the same strategy than in [B.84a] for eigenvalues \ $\not= 1$ \ and [B.84 b] for the eigenvalue 1.$\hfill \blacksquare$

\newpage

\section*{References.}

\bigskip

\bigskip

\begin{enumerate}

\item{[B.84 a)]} Barlet, D. \textit{Contribution effective de la monodromie aux d\'eveloppements asymptotiques,} Ann. Sc. Ec. Norm. Sup. 17 (1984), p.293-315.

\item{[B.84.b]} Barlet, D.  \textit{Contribution du cup-produit de la fibre de Milnor aux p\^oles de \ $\vert f \vert^{2\lambda}$ }, Ann. Inst. Fourier (Grenoble), 34 (1984) p.75-107.

 \item{[B.93]} Barlet, D. \textit{Theory of (a,b)-Modules I,} Complex Analysis and Geometry, Plenum Press New York, (1993), p.1-43.
 
 \item{[B.95]} Barlet, D. \textit{Theorie des (a,b)-modules II. Extensions}, Complex Analysis and Geometry, Pitman Research Notes in Math. Series 366,  (Trento 1995), p.19-59, Longman (1997).

 \item{[B.04]} Barlet, D.  \textit{Sur certaines singularit\'es non isol\'ees d'hypersurfaces I}, preprint de l'Institut E. Cartan (Nancy) 2004/n$^{\circ}$03, 47 pages. A second version (shorter) will appear in Bull. Soc. Math. France.
 
\item{[B.05]} Barlet, D.  \textit{Sur certaines singularit\'es non isol\'ees d'hypersurfaces II}, preprint de l'Institut E. Cartan (Nancy) 2005/n$^{\circ}$42, 47 pages.
 
 \item{[Be.01]} Belgrade, R. \textit{Dualit\'e et Spectre des (a,b)-modules,} Journal of Algebra 245, (2001), p.193-224.
 
 \item{[Br.70]} E. Brieskorn : {\it Die Monodromie der isolierten Singularit\"aten von Hyperfl\"achen.} Manuscripta Math. 2 (1970), p. 103-161.
 
 \item{[K.76]} Kashiwara, M. \textit{b-Function and Holonomic Systems, Rationality of Roots of b-Functions}, Invent. Math. 38 (1976) p.33-53.

\item{[M.74]} Malgrange, B. \textit{Int\'egrales asymptotiques et monodromie,} Ann. Sc. Ec. Norm. Sup. , t.7, (1974), p.405-430.

\end{enumerate}

 \bigskip

\noindent Daniel Barlet, 

\noindent Universit\'e Henri Poincar\'e (Nancy I ) et Institut Universitaire de France,

\noindent Institut E.Cartan  UHP/CNRS/INRIA, UMR 7502 ,

\noindent Facult\'e des Sciences et Techniques, B.P. 239

\noindent 54506 Vandoeuvre-les-Nancy Cedex , France.

\noindent e-mail :  barlet@iecn.u-nancy.fr

 \end{document}